\documentclass[11pt,draft]{article}

\usepackage[T1]{fontenc}
\usepackage{amsfonts}

\newtheorem{thm}{theorem}[section]
\newtheorem{theorem}[thm]{Theorem}
\newtheorem{proposition}[thm]{Proposition}
\newtheorem{lemma}[thm]{Lemma}

\newtheorem{example}[thm]{Example}
\newtheorem{definition}[thm]{Definition}

\begin{document}

\title{Primeness property for central polynomials of verbally prime P.I. algebras}

\author{Diogo Diniz Pereira da Silva e Silva \footnote{Supported by CNPq}}
\date{{\small \textit{Unidade Acadêmica de Matemática e Estatística\\
Universidade Federal de Campina Grande\\
Cx. P. 10.044, 58429-970, Campina Grande, PB, Brazil}}\\
E-mail: diogo@dme.ufcg.edu.br}
\maketitle

\textbf{Keywords}: verbally prime algebras, central polynomials, regular algebras

\textbf{MSC}: 16R10 16R99 15A75

\maketitle

\begin{abstract}
Let $f$ and $g$ be two noncommutative polynomials in disjoint sets of variables. An algebra $A$ is verbally prime if whenever $f\cdot g$ is an identity for $A$ then either $f$ or $g$ is also an identity. As an analogue of this property Regev proved that the verbally prime algebra $M_k(F)$ of $k\times k$ matrices over an infinite field $F$ has the following primeness property for central polynomials:  whenever the product $f\cdot g$ is a central polynomial for $M_k(F)$ then both $f$ and $g$ are central polynomials. In this paper we prove that over a field of characteristic zero Regev' s result holds for the verbally prime algebras $M_k(E)$ and $M_{k,k}(E)$, where $E$ is the infinite dimensional Grassmann algebra.
\end{abstract}

\section{Introduction}

The notion of verbally prime algebras was introduced by Kemer (see \cite{Kemer}) and these algebras are of fundamental importance in PI-theory. An algebra $A$ is verbally prime if whenever the product of two polynomials $f$ and $g$ in distinct variables is an identity for $A$ then either $f$ or $g$ is also an identity. In the case of characteristic zero the verbally prime algebras are (up to p.i. equivalence) $M_n(K)$, the algebra of $n\times n$ matrices over the field $K$, the algebra $M_n(E)$ of $n\times n$ matrices over the Grassmann algebra of an infinite dimensional vector space and its subalgebras $M_{p,q}(E)$, where 
$p+q=n$. 

The existence of central polynomials for verbally prime algebras was proved by Kemer \cite{Kemer2} and earlier for matrix algebras independently by Formanek \cite{Formanek} and Razmyslov 
\cite{Razmyslov}. As a consequence of this some important results were obtained using central polynomials for matrix algebras (see \cite[Chapter 4]{DrenskyFormanek}). The existence of central polynomials for verbally prime algebras was used by A. Belov in \cite{Belov} to prove that each associative PI-algebra does not coincide with its commutant. 

Amitsur proved in \cite{Amitsur} that $M_k(F)$ is verbally prime. In fact a stronger statement was proved: if the product of any two polynomials $f\cdot g$ is an identity for $M_k(F)$ then either $f$ or $g$ is. An analogous property for the central polynomials of $A=M_k(F)$ was proved by Regev in \cite{Regev}: if the product $f(x_1,\dots, x_r)\cdot g(x_{r+1},\dots, x_s)$ is a central polynomial for $A$ then both $f$ and $g$ are central polynomials. 

In this paper we prove that Regev's result holds for the verbally prime algebras $M_k(E)$ and $M_{k,k}(E)$. We use the regular decompositions of the algebras $M_k(E)$ and 
$M_{k,k}(E)$. Regular algebras were defined in \cite{RegevSeeman} in the study of multilinear identities of $\mathbb{Z}_2$-graded tensor products of p.i. algebras. 

\section{Preliminaries}

\subsection{Identities and Central Polynomials}
Let $F$ be a field of characteristic zero. In this article all algebras, vector spaces and tensor products are over $F$. Denote by $F\langle X \rangle$ the free associative algebra on the free generating set $X=\{x_1,x_2,\dots\}$. The degree of a polynomial 
$f$ in the variable $x_i$ is denoted by $deg_{x_i}(f)$. The polynomial $f(x_1,\dots, x_n)$ is \textit{multilinear} if in every monomial in $f$ the degree of $x_i$ is $1$ for $1\leq i \leq n$. An ideal $I$ of $F\langle X \rangle$ is a \textit{$T$-ideal} if it is invariant under all endomorphisms.

A polynomial $f(x_1,\dots, x_n)$ is a \textit{polynomial identity} for an algebra $A$ if 
$f(a_1,\dots, a_n)=0$ for every $a_1,\dots, a_n \in A$. The center of $A$ is denoted by $Z(A)$, and if $f(a_1,\dots, a_n) \in Z(A)$ for every $a_1,\dots, a_n$ then $f$ is a \textit{central polynomial} for $A$. Moreover if $f$ is not an identity then it is a \textit{proper central polynomial} for $A$. The set $T(A)$ of all polynomial identities for $A$ is a $T$-ideal of $A$. In order to simplify the statments of the results in this article we make the following:

\begin{definition}
If an algebra $A$ has proper central polynomials for $A$ and satisfies the following property: whenever 
$f(x_1,\dots, x_r)\cdot g(x_{r+1},\dots, x_s)$ is a proper central polynomial for $A$ then both $f$ and $g$ are proper central polynomials for $A$, we say that $A$ has the primeness property on central polynomials.
\end{definition}

\subsection{Graded Algebras and Verbally prime algebras}

Let $G$ be an arbitrary group. The decomposition $A=\oplus_{g\in G}A_g$ is a \textit{$G$-grading}  if $A_gA_h\subset A_{gh}$. A non-zero element in $\cup_{g\in G} A_g$ is said to be homogeneous of degree $g$ if $a \in A_g$ and its degree is denoted by $deg(a)$. 

Now we describe the verbally prime algebras. Let $V$ be a vector space with basis $\{e_1,e_2,\dots\}$. We denote by $E$ the Grassmann algebra of $V$ and recall that this algebra has the presentation: \[E=\langle 1,e_1, e_2, \dots | e_ie_j=-e_je_i,\mbox{ for all }i,j \geq 1 \rangle.\] It is well known that $B=\{1, e_{i_1}\dots e_{i_k}| 1\leq i_1<\dots < i_k\}$ is a basis of $E$. Moreover if $E_0$ is the subspace generated by $1$ and the monomials in $B$ with even $k$ and $E_1$ is the subspace generated by the monomials in $B$ with odd $k$ then 
$E=E_0\oplus E_1$ is a $\mathbb{Z}_2$-grading for $E$. Moreover $Z(E)=E_0$ and the non-zero product of two homogeneous elements is central if and only if they have the same degree. Kemer proved that every nontrivial verbally prime algebra is p.i. equivalent to one of the algebras $M_k(F)$, $M_{k}(E)$ or to the subalgebra $M_{p,q}(E)$ of $M_{p+q}(E)$ that we describe in the next paragraph.

Let $(M_{p+q}(F))_0$ be the subspace of matrices $\left(\begin{array}{cc}
	a&0\\
	0&b
\end{array}
 \right)$, where $a \in M_p(F)$ and $b \in M_{q}(F)$ and  $(M_{p+q}(F))_1$ the subspace of matrices
 $\left(\begin{array}{cc}
	0&c\\
	d&0
\end{array}
 \right)$, where $c \in M_{p\times q}(F)$ and $d \in M_{q\times p}(F)$. This is a $2$-grading of $M_{p+q}(F)$ and $(M_{p+q})_0\otimes E_0\oplus (M_{p+q})_1\otimes E_1$ is a subalgebra of $M_{p+q}(F)\otimes E$. This tensor product is isomorphic to $M_{p+q}(E)$ and we denote by $M_{p,q}(E)$ the corresponding subalgebra.

\subsection{Regular Algebras}

The following definition of regular algebras was given in \cite[Definition 2.3]{RegevSeeman}.

\begin{definition}
Let $A$ be an algebra and consider the decomposition $A=A_1\oplus \dots \oplus A_r$, where the $A_i$ are vector spaces. Such a decomposition is called "regular" if it satisfies:
\begin{enumerate}
\item[(P1)] Given $n$ and an ordered sequence $i_1,\dots, i_n$ ($1\leq i_j \leq r$ for all $i_j$), there exist $x_1,\dots, x_n$ where each $x_j \in A_{i_j}$, and $x_1\dots x_n \neq 0$.
\item[(P2)] For all $x \in A_i$, $y \in A_j$ $(1\leq i, j \leq r)$, $xy = \epsilon_{i,j}^{A} yx$, where $0\neq \epsilon_{i,j}^{A}\in F$ depends only on $i$ and $j$. Thus for homogeneous $x,y \in A$, 
\begin{equation}\label{Commutation}
xy=\epsilon_{deg(x),deg(y)}^{A}yx.
\end{equation}
\end{enumerate}
\end{definition}



Denote by $M^{A}$ the $r\times r$ matrix \[M^{A}=(\epsilon_{i,j}^{A})_{1\leq i,j \leq r}.\] In this case we say that the algebra $A$ is \textit{regular} and that the matrix $M^{A}$ is the \textit{matrix of commutation relations} of the regular decomposition. 

\begin{example}\cite[Example 1.2]{BahturinRegev}
The algebra $M_k(F)$ is regular. Let $\xi$ be a primitive $k$-th root of $1$ and assume that $\xi \in F$. Fix the matrices $X_a=diag(\xi^{k-1},\xi^{k-2},\dots, 1)$ and 
$X_b=e_{n,1}+\sum_{i=1}^{k-1}e_{i,i+1}$, where $e_{i,j}$ are the matrix units. Then \[M_{k}(F)=\bigoplus_{i,j=0}^{k-1}A_{i,j},\mbox{ where }A_{i,j}=F\cdot (X_a^{i}X_b^j),\] is a regular decomposition.
\end{example}

\begin{example}\cite[Example 2.3]{BahturinRegev}\label{Example}
The algebra $M_{1,1}(E)$ is regular. Let \[I=\left(\begin{array}{cc}
	1&0\\
	0&1
\end{array}
 \right), X_a=\left(\begin{array}{cc}
	-1&0\\
	0&1
\end{array}
 \right), X_b\left(\begin{array}{cc}
	0&1\\
	1&0
\end{array}
 \right)\mbox{ and }X_aX_b=\left(\begin{array}{cc}
	0&-1\\
	1&0
\end{array}
 \right).\] The decomposition \[M_{1,1}(E)=E_0I\oplus E_0X_a \oplus E_1X_b \oplus E_1X_aX_b,\] is regular with matrix \[M^{M_{1,1}(E)}=\left(\begin{array}{cccc}
	1&1&1&1\\
	1&1&-1&-1\\
	1&-1&-1&1\\
	1&-1&1&-1\\
\end{array}
 \right).\]
This regular decomposition is a $\mathbb{Z}_2\times \mathbb{Z}_2$-grading and the center of $M_{1,1}(E)$ is the $(0,0)$ component. The components $(0,1)$, $(1,0)$ and $(1,1)$ are  
$E_0X_a$, $E_1X_b$ and $ E_1X_aX_b$ respectively.
\end{example}

The tensor product of regular algebras has a natural regular decomposition, more precisely: 

\begin{proposition}\cite[Proposition 2.1]{BahturinRegev}
Let \[A=\bigoplus_{i=1}^{r} A_{i}\mbox{ be regular with }M^{A}=(\theta_{i,j})_{1\leq i,j \leq r.}\]
Similarly let \[B=\bigoplus_{i=1}^{s} B_{i}\mbox{ be regular with }M^{B}=(\eta_{k,l})_{1\leq k,l \leq s.}\]
The ordinary tensor product $A\otimes B$ is regular for the decomposition \[A\otimes B= \bigoplus_{i=1}^{r}\bigoplus_{k=1}^{s}A_i\otimes B_k.\]
It has the matrix of commutation relations $M^{A\otimes B}=M^{A}\otimes M^{B}$.
\end{proposition}
In the next setion we prove that if $A$ and $B$ are verbally prime regular algebras with the primeness property on central polynomials then $A\otimes B$ also has the property. 
The next remarks are used in the proof of this result. Let $a=(a_1,\dots, a_n)$ be a sequence of non-zero elements in $A$, it follows from the commutation relations (\ref{Commutation}) that for every permutation $\sigma \in S_n$ there exists $0\neq \epsilon_{\sigma}^{A} \in F$ that depends on $a$ and $\sigma$ only such that \[a_{\sigma(1)}\dots a_{\sigma(n)}= \epsilon_{\sigma}^{A} a_1\dots a_n.\] For example, the  
$2$-grading of the Grassmann algebra $E$ is a regular decomposition and if $a_i =e_i$ for every $i$ then $\epsilon_{\sigma}^{A}$ is just the sign of the permutation 
(see \cite[Example 1.1.8]{GiambrunoZaicev}). Let $A=E\otimes E$ and let $C\cup D$ be a partition of $\{1,2,\dots, n\}$. Define $a_i=e_{i}\otimes 1$ if $i \in C$ and $a_i=1\otimes e_i$ if $i \in D$. Then $\epsilon_{\sigma}^{A}$ is just the coloured sign of $\sigma$ (see \cite[Corollary20]{PlamenAzevedo}).

\begin{definition}
Let $f(x_1,\dots, x_n)=\sum_{\sigma \in S_n}\alpha_{\sigma}x_{\sigma(1)}\dots x_{\sigma(n)}$ be a multilinear polynomial in $F\langle X \rangle$ and $a=(a_1,\dots, a_n)$ be a sequence of elements in the regular algebra $A$. We denote by $f_a$ the polynomial 
\[f_a(x_1,\dots,x_n)=\sum_{\sigma \in S_n}\epsilon_{\sigma}^{A}\alpha_{\sigma}x_{\sigma(1)}\dots x_{\sigma(n)}.\]
\end{definition}

The following proposition is an easy consequence of the previous discussion and is used in the proof of Theorem \ref{TR}.

\begin{proposition}\label{propp}\label{PP}
Let $A$ and $B$ be regular algebras and $f(x_1,\dots, x_n)=\sum \lambda_{\sigma} x_{\sigma(1)}\dots x_{\sigma(n)}$ be a multilinear polynomial. For every sequence $a=(a_1,\dots, a_n)$ of elements in $A$ and every sequence $b=(b_1,\dots, b_n)$of elements in $B$ we have 
\[f(a_1\otimes b_1,\dots, a_n\otimes b_n)=(a_1\dots a_n)\otimes f_a(b_1,\dots, b_n)=f_b(a_1\dots a_n)\otimes b_1,\dots, b_n.\]
\end{proposition}
\textit{Proof}
 Note that 
\begin{eqnarray*}
f(a_1\otimes b_1,\dots, a_n\otimes b_n)&=&\sum \lambda_{\sigma} (a_{\sigma(1)}\dots a_{\sigma(n)}\otimes b_{\sigma(1)}\dots b_{\sigma(n)})\\
                                       &=&\sum \lambda_{\sigma}\epsilon_{\sigma}^{A} (a_{1}\dots a_{n}\otimes b_{\sigma(1)}\dots b_{\sigma(n)})\\
																			&=&a_{1}\dots a_{n}\otimes f_a(b_1,\dots, b_n).														
\end{eqnarray*}
The proof of the second equality is analogous.
\hfill $\Box$

\section{Results}
In this section we prove the main results. 

\begin{proposition}\label{PROP}
Let $A$ be an algebra with proper central polynomials. Supose that whenever $f,g$ are  multilinear polynomials in distinct variables such that 
$f(x_1,\dots, x_r)\cdot g(x_{r+1},\dots, x_s)$ is a proper central polynomial for $A$ then both $f$ and $g$ are proper central polynomials for $A$. Then the algebra $A$ has the primeness property on central polynomials.
\end{proposition}

\textit{Proof}
 Let $f,g$ be  arbitrary polynomials in distinct variables such that $f(x_1,\dots, x_r)\cdot g(x_{r+1},\dots, x_s)$ is a proper central polynomial for $A$. Write 
\[f(x_1,\dots, x_r)=\sum f_i(x_1,\dots, x_r) \mbox{ and } g(x_{r+1},\dots, x_s)=\sum g_j(x_{r+1},\dots, x_s)\] as a sum of its multihomogeneous components. The  products $f_i\cdot g_j$ are the multihomogeneous components of $f\cdot g$ and therefore $f_i\cdot g_j$ is a central polynomial for $A$ for every $i,j$. We apply the multilinearization process to $f_i\cdot g_j$ and since the polynomials are in distinct variables we obtain the multilinear polynomial $f_i^{\prime}\cdot g_j^{\prime}$ where $f_i^{\prime}$, $ g_j^{\prime}$ are the complete multilinearizations of $f_i$ and $g_j$ respectively. Note that  $f_i^{\prime}\cdot g_j^{\prime}$  is a proper central polynomial for $A$ and therefore both $f_i^{\prime}$ and $ g_j^{\prime}$ are central polynomials. Since the field $F$ is of characteristic zero this implies that both $f_i$ and $g_j$ are central polynomials. Hence we conclude that $f$ and $g$ are central polynomials for $A$.
\hfill $\Box$

\begin{theorem}\label{TR}
Let $A$ and $B$ be verbally prime algebras and assume that these are regular algebras. If $A$ and $B$ have the primeness proprety on central polynomials then $A\otimes B$ has the same property.
\end{theorem}
\textit{Proof}
Let $f(x_1,\dots, x_r)$, $g(x_{r+1},\dots, x_s)$ be polynomials in disjoint sets of variables such that $f\cdot g$ is a central polynomial for 
$C=A\otimes B$. It follows from Proposition \ref{PROP} that we may assume that both polynomials are multilinear. We prove that $f$ is a proper central polynomial for $C$, the proof for $g$ is analogous. Since $g$ is a multilinear polynomial and is not an identity there exists $c_i^{\prime}=(a_{i}^{\prime}\otimes b_{i}^{\prime})$ in $C$, $r+1\leq i \leq s$, such that $g(c_{r+1}^{\prime},\dots, c_s^{\prime})\neq 0$. Clearly this implies that 
$g_{a^{\prime}}(b_{r+1}^{\prime},\dots, b_s^{\prime})\neq 0$, where $a^{\prime}=(a_{r+1}^{\prime},\dots, a_s^{\prime})$. Therefore $g_{a^{\prime}}(x_{r+1},\dots, x_s)$ is not an identity for $B$. We claim that for an arbitrary $r$-tuple $a=(a_1,\dots,a_r)$ of elements in $A$ the corresponding polynomial $f_a(x_1,\dots, x_r)$ is a central polynomial for $B$. Since $B$ has the primeness property and it is verbally prime to prove this claim it is enough to prove that the product $f_a(x_1,\dots, x_r)\cdot g_{a^{\prime}}(x_{r+1},\dots, x_s)$ is a central polynomial for $B$. Given 
$b_i$, $1\leq i \leq s$, we get that 
\begin{eqnarray*}
&&f(a_1\otimes b_1,\dots, a_r\otimes b_r)\cdot g(a_{r+1}^{\prime}\otimes b_{r+1},\dots, a_{r+1}^{\prime}\otimes b_{r+1})
\\&=&a_1\dots a_ra_{r+1}^{\prime}\dots a_{s}^{\prime}\otimes f_a(b_1,\dots,b_r)g_{a^{\prime}}(b_{r+1},\dots,b_{s}),
\end{eqnarray*}
lies in $Z(C)$. Therefore $f_a(b_1,\dots,b_r)g_{a^{\prime}}(b_{r+1},\dots,b_{s})$ lies in $Z(B)$ and the claim is proved. Analogously for an arbitrary $r$-tuple $b=(b_1,\dots,b_r)$ the corresponding polynomial $f_b(x_1,\dots, x_r)$ is a central polynomial for $A$. To prove that $f$ is a central polynomial for $C$ it is enough to consider substitutions of the variables by elementary tensors. Let $c_i=(a_{i}\otimes b_{i})$, 
$1\leq i \leq r$, where $a_i \in A$ and $b_i \in B$. It follows from Proposition \ref{PP} that 
\[f(a_1\otimes b_1,\dots, a_n\otimes b_n)=(a_1\dots a_n)\otimes f_a(b_1,\dots, b_n)=f_b(a_1\dots a_n)\otimes b_1,\dots, b_n.\] Since $f_a(b_1,\dots, b_n)$ lies in $Z(B)$ and 
$f_b(a_1\dots a_n)$ lies in $Z(A)$ we conclude that $f(a_1\otimes b_1,\dots, a_n\otimes b_n)$ lies in $Z(C)$. Clearly $f$ is not an identity and therefore is a proper central polynomial for $C$.
\hfill $\Box$

Our objective now is to prove that the Grassman algebra has the primeness property on central polynomials. The proof uses the well known description of the polynomial identities of $E$ (see for example \cite{DrenskyFormanek}).

\begin{lemma}\cite[Lemma 1.4.3]{DrenskyFormanek}
Modulo the polynomial identity \\$[x_1,x_2,x_3]=0$, every multilinear polynomial can be presented in the form 
\begin{equation}\label{EE}
f(x_1,\dots, x_m)=\sum \beta_{i,j}x_{i_1}\dots x_{i_{m-2p}}[x_{j_1},x_{j_2}]\dots[x_{j_{2p-1}},x_{j_{2p}}], \beta_{i,j}\in K,
\end{equation}
where the sum runs on all permutations $[i_1,\dots, i_{2m-p},j_1,j_2,\dots,j_{2p-1},j_{2p}]$ of $1,\dots, m$ such that $i_1<\dots < i_{m-2p}$, $j_1<j_2<\dots<j_{2p-1}<j_{2p}$.
\end{lemma}

The proof of the next lemma is based on the proof of \cite[Theorem 1.4.4]{DrenskyFormanek}.

\begin{lemma}\label{EEE}
If $f(x_1,\dots, x_m)$ is a multilinear polynomial that is not an identity for $E$ then there existis $a_1,\dots, a_n$ in $E$ such that $f(a_1,\dots, a_n)$ is a non-zero element in the center of $E$.
\end{lemma}

\textit{Proof}
We may assume that $f(x_1,\dots, x_m)$ is written as in (\ref{EE}). We choose the minimal $p$ with $\beta_{i,j}\neq 0$ and fix one of the corresponding permutations $[i_1,\dots, i_{2m-p},j_1,j_2,\dots,j_{2p-1},j_{2p}]$. We consider the following products of two generators of the Grassmann algebra $a_{i_1}=e_1e_2, \dots, a_{i_{m-2p}}=e_{2m-4p-1}e_{2m-4p-1}$ and the following generators 
$a_{j_1}=e_{2m-4p+1}$, \\$a_{j_2}=e_{2m-4p+2}$,$\dots$, $a_{j_{2p-1}}=e_{2m-2p-1}$,$a_{j_{2p}}=e_{2m-2p-1}$. Note that 
\begin{eqnarray*}
&&a_{i_1}\dots a_{i_{m-2p}}[a_{j_1},a_{j_2}]\dots [a_{j_{2p-1}},a_{j_{2p}}]
\\&=&2^{p}e_1e_2\dots e_{2m-4p-1}e_{2m-4p}e_{2m-4p+1}e_{2m-4p+2}\dots e_{2m-4p+1}e_{2m-4p+2}\dots e_{2m-2p-2}e_{2m-2p}
\end{eqnarray*}
is a non-zero element in the center of $E$. Moreover all products $e_ae_b$ lie in the center, and since $p$ is minimal we conclude that all other summands
$a_{i_1^{\prime}}\dots a_{i_{m-2p}^{\prime}}[a_{j_1^{\prime}},a_{j_2^{\prime}}]\dots [a_{j_{2p-1}^{\prime}},a_{j_{2p}^{\prime}}]$ are zero and the result of this substitution is a non-zero element in the center of $E$.
\hfill $\Box$

\begin{proposition}\label{LE}
The infinite dimensional Grassmann algebra $E$ has the primeness property on central polynomials.
\end{proposition}

\textit{Proof}
Let $f(x_1,\dots, x_r)$, $g(x_{r+1},\dots, x_s)$ be polynomials in disjoint sets of variables such that $f\cdot g$ is a proper central polynomial for $E$. We may assume that both polynomials are multilinear. It follows from Lemma \ref{EEE} that there exists $a_{r+1},\dots, a_s$ in $E$ such that $b=g(a_{r+1},\dots, a_s)\neq 0$ has degree $0$ in the 
$\mathbb{Z}_2$-grading of $E$. For every $n$ the subalgebra $E^{\prime}$ of $E$ generated by the $e_j$ with $j\geq n$ is isomorphic to $E$. Let $n$ be large enough such that if 
$a^{\prime}$ is a nonzero element of $E^{\prime}$ then $a^{\prime}\cdot b\neq 0$. We prove that $f$ is a central polynomial for $E^{\prime}$ and it follows that $f$ is a central polynomial for $E$. Given $a_1,\dots, a_r$ in $E^{\prime}$ the product 
$f(a_1,\dots, a_r)\cdot b$ lies in $Z(E)=E_0$. If $f(a_1,\dots,a_r)\neq 0$ then  $f(a_1,\dots, a_r)\cdot b\neq 0 $ and we conclude that $f(a_1,\dots, a_r)$ has degree zero. Since 
$Z(E^{\prime})=E_0\cap E^{\prime}$ the polynomial $f$ is a central polynomial for $E^{\prime}$. The proof that $g$ is a central polynomial is analogous.
\hfill $\Box$



%

We recall the theorem of Regev on the central polynomials of $M_k(F)$.

\begin{theorem}\cite[Regev]{Regev}\label{TRegev}
Lef $F$ be an infinite field. The algebra $M_k(F)$ has the primeness property on central polynomials
\end{theorem}

Next we prove our main result.

\begin{theorem}
The verbally prime regular algebras $M_k(E)$ and $M_{k,k}(E)$ have the primeness property on central polynomials.
\end{theorem}

\textit{Proof}
As a consequence of Proposition \ref{LE} and Theorem \ref{TR} we conclude that the tensor product $E\otimes E$ has the primeness property. The algebra $M_{1,1}(E)$ satisfies the same polynomials identities as $E\otimes E$ and therefore has the same central polynomials. This implies that the algebra $M_{1,1}(E)$ also has the primeness property. Finally the result follows from the isomorphisms $M_k(E) \cong M_k(F)\otimes E$ and $M_{k,k}(E)\cong M_k(E)\otimes M_{1,1}(E)$ together with Theorem \ref{TR} and Theorem \ref{TRegev}.
\hfill $\Box$

We remark finally that it is interesting to study the primeness property on the central polynomials of the verbally prime algebras $M_{p,q}(E)$ over a field of characteristic zero and to consider the algebras $M_n(E)$ and $M_{p,q}(E)$ over infinite fields.

\begin{flushleft}
\textbf{Acknowledgements}
\end{flushleft}
We thank CNPq for the finantial support.

\end{document}